\documentclass[amssymb,twoside,10pt,amscd,reqno]{amsart}

\usepackage{latexsym,amsfonts}
\usepackage{amsmath,amsthm}          
\usepackage{amssymb}                 
\usepackage{amscd}                   
\usepackage{euscript}                

\def\segp{\mathrel{\widehat\otimes}}
\def\oM{\overline{M}}
\def\oQ{\overline{Q}}

\def\bP{\Bbb P}

\def\bV{\Bbb V}

\def\Aut{\operatorname{Aut}}

\def\Sym{\operatorname{Sym}}

\def\Hilb{\operatorname{Hilb}}

\def\Tor{\operatorname{Tor}}

\def\sym{\Sym}

\def\cO{\Cal O}
\def\cS{\Cal S}
\def\cD{\Cal D}

\def\codim{\operatorname{codim}}
\def\deg{\operatorname{deg}}

\def\bZ{\Bbb Z}
\def\bC{\Bbb C}
\def\bQ{\Bbb Q}

\def\cD{\Cal D}
\def\cS{\Cal S}

\def\cF{\Cal F}

\def\cM{\Cal M}

\def\Ker{\operatorname{Ker}}
\def\Im{\operatorname{Im
}}

\def\Cal{\mathcal}

\newtheoremstyle{mystyle}{}{}{\itshape}{}{\scshape}{.}{ }{}
\theoremstyle{mystyle}

\swapnumbers

\newtheorem{Theorem}{Theorem}[section]

\newtheorem{Proposition}[Theorem]{Proposition}
\newtheorem{Lemma}[Theorem]{Lemma}
\newtheorem{Corollary}[Theorem]{Corollary}

\newtheoremstyle{myreview}{}{}{}{}{\scshape}{.}{ }{}
\theoremstyle{myreview}

\newtheorem{Definition}[Theorem]{Definition}
\newtheorem{Example}[Theorem]{Example}

\newtheorem{Remark}[Theorem]{Remark}

\newtheorem{Notation}[Theorem]{Notation}

\newcounter{et}[Theorem]
\def\cooltag{\tag{\arabic{section}.\arabic{Theorem}.\arabic{et}}\addtocounter{et}{1}}

\def\boplus{\bigoplus\limits}
\def\cL{\Cal L}
\def\cG{\Cal G}

\begin{document}

\title{Equations for $\oM_{0,n}$}
\author{Sean Keel and Jenia Tevelev}
\address{Department of Mathematics, University of Texas at Austin,
Austin, Texas, 78712}
\email{keel@math.utexas.edu and tevelev@math.utexas.edu}

\begin{abstract} We show that the log canonical bundle, $\kappa$, of $\oM_{0,n}$
is very ample, show the homogeneous coordinate ring is Koszul, and give a nice
set of rank $4$ quadratic generators for the homogeneous ideal: The embedding
is equivariant for the symmetric group, and the image lies on many
Segre embedded copies of
$\bP^1 \times \dots \times \bP^{n-3}$, permuted by the symmetric group.
The homogeneous ideal of $\oM_{0,n}$
is the sum of the homogeneous ideals of these Segre embeddings.
\end{abstract}

\maketitle

\section{Introduction and Statement of Results}

Let $\oM_S = \oM_{0,n}$ be the moduli space of stable $n$-pointed
curves, with marked points labeled by the elements of the finite
set $S$, with cardinality $|S| = n$,
over an algebraically closed field $k$ of characteristic zero.
Our goal is to study the equations of $\oM_S$ in its most natural embedding.
For a finite set $T$, let $W_T$ be the standard irreducible
representation of the symmetric group $\Aut(T)$, i.e. $T$-tuples
of integers that sum to zero.

\begin{Theorem}\label{maintheorem}
The log canonical line bundle
$$
\kappa := \cO(K_{\oM_S} + B)
$$
is very ample on $\oM_S$, where
$B:=\oM_S \setminus M_S$ is the boundary.
A flag of subsets
$$
S_3 \subset S_4 \subset \dots \subset S_n = S
$$
of $S$ with $|S_k| = k$ canonically induces an identification
$$
H^0(\oM_S,\kappa) = W_{S_3}\otimes W_{S_5}\otimes\ldots\otimes W_{S_{n-1}}
$$
and $\oM_S \subset \bP(H^0(\oM_S,\kappa)^*)$
factors through a Segre embedding
$$
\bP(W_{S_3}^*) \times \bP(W_{S_4}^*)\times\ldots \times \bP(W_{S_{n-1}}^*)
\subset \bP(H^0(\oM_S,\kappa)^*).
$$
The homogeneous ideal of $\oM_S$ is the sum of the homogeneous
ideals for all
Segre embeddings over all flags of subsets, i.e. we have
the ideal theoretic equality
$$
\oM_S = \bigcap_{S_3 \subset S_4 \dots \subset S}
\bP(W_{S_3}^*) \times \bP(W_{S_4}^*)\times\ldots \times \bP(W_{S_{n-1}}^*)
\subset \bP(H^0(\oM_S,\kappa)^*).
$$
\end{Theorem}

The coordinate ring $R:= \boplus_{n \geq 0} H^0(\oM_S,\kappa^{\otimes n})$
has very nice properties:

\begin{Theorem}\label{propertiesofkapparing}
$R$ is Koszul, i.e. the trivial $R$-module $k$ has a resolution
of form
$$
\dots \to R[-k]^{a_k} \to R[-(k-1)]^{a_{k-1}} \to \dots \to R[-1]^{a_1}
\to R \to k \to 0.
$$
$\oM_S$ is projectively normal, i.e.~the natural map
$$
\Sym^\bullet(R_1) \rightarrow R
$$
is surjective. Its kernel is generated by quadrics of rank at most $4$.
\end{Theorem}

The embedding
$$
\oM_S \subset \bP:=\bP^1 \times \dots \times \bP^{n-3}
$$
has nice properties as well.
Let
$$
\cL := \cO_{\bP^1}(1) \boxtimes \dots \boxtimes \cO_{\bP^{n-3}}(1)
$$
and let $B$ be the coordinate ring
$$
B:= \boplus_{n \geq 0} H^0(\bP, \cL^{\otimes n}).
$$

\begin{Theorem}\label{npforms}
$R$ is the homogeneous coordinate ring of $\oM_S \subset \bP$. The
embedding satisfies the analog of Green--Lazarsfeld's property $N_p$
for all $p$ i.e. the minimal resolution of $R$ over $B$ is of form
$$
\dots \to B[-k]^{a_k} \to \dots \to B[-2]^{a_2} \to B \to R \to 0.
$$
\end{Theorem}
\begin{Remark} The analogous statement for
$$
\oM_S \subset \bP(H^0(\kappa)^*)
$$
fails. For example $\oM_{0,5} \subset \bP^5$ fails to satisfy $N_3$. In
fact the only non-degenerate irreducible subvarieties $X$ of projective space that
satisfy $N_p$ for all $p$ are varieties of minimal
degree $\deg X=1+\codim X$, i.e.~quadric hypersurfaces,
rational normal scrolls, or a cone over the Veronese surface in $\bP^5$,
see~\cite{EH,EG}.
\end{Remark}

We note that the compactification $M_S \subset \oM_S$ is canonical,
indeed the so called log canonical compactification of (the log
minimal variety) $M_S$, see \cite{KeelTevelev}, and the coordinate ring
$R$ is the  log canonical ring of $M_S$.

The above results allow us
to construct $\oM_S$ from $M_S$ in a purely topological and
combinatorial way:
Dropping the last point gives a tautological
fibration
$$
\pi:M_{0,n} \rightarrow M_{0,n-1}
$$
with fibre $\bP^1\setminus \{p_1,\ldots,p_{n-1}\}$
for the distinct points $p_1,\ldots,p_{n-1}$.

Arnold \cite{Ar}
observed that the Leray spectral sequence degenerates at the $E^2$ term
inducing a canonical
isomorphism
$$
H^{n-3}(M_{0,n},\bZ) = H^{n-4}(M_{0,n-1},\bZ)
\otimes H^1(\bP^1 \setminus \{p_1,\ldots,p_{n-1}\},\bZ)
$$
See Theorem~\ref{old2.5}.
By induction we have
$$
H^{n-3}(M_{0,n},\bZ) = \bigotimes_{i = 3}^{{n-1}}
H^1(\bP^1 \setminus \{p_1,\dots,p_i\},\bZ).
$$
Note that $H^1(\bP^1 \setminus\{p_1,\dots,p_i\},\bZ)$ is
naturally identified with $W_i$, the standard irreducible representation
of the symmetric group~$\Aut\{p_1,\ldots,p_i\}$.
We have in particular a tautological identification
$$
H_{n-3}(M_{0,n},\bZ) = W_3^*\otimes\ldots\otimes W_{n-1}^*.
$$
The symmetric group $S_n$ acts naturally on $M_{0,n}$,
and so it acts on the left hand side of
the above equality.
This is the so-called Whitehouse module, introduced by Kontsevich \cite{K}
in the context of (cyclic) Lie operads.
This action does not preserve the tensor product,
but by~\ref{maintheorem} (see Theorem \ref{old2.5}) we have

\begin{Corollary}
Let $X$ be the set of totally decomposable elements (i.e. elements of form
$x_3 \otimes x_2\otimes \dots \otimes x_{n-1}$ that remain
decomposable under the symmetric group).
$X$ is closed under scalar multiplication
and so defines a subset
$$\bP(X)\subset \bP(H_{n-3}(M_{0,n},\bC)).$$
This subset is canonically identified with
$\oM_{0,n}$, embedded by the full
log canonical
series
$
|K_{\oM_{0,n}} + \partial \oM_{0,n}|
$
of $M_{0,n}$.
\end{Corollary}

Thus in particular we see that the compactification
$M_{0,n} \subset \oM_{0,n}$
can be canonically recovered from the homology of
the spaces $M_{0,k}$ together with their symmetric
group action, i.e.~from the cyclic Lie operad \cite{Getzler}.

The paper is organized as follows: \ref{segrequadrics}  contains the
proof of Theorem \ref{maintheorem}. \ref{kappa}~contains the proof
of Theorems \ref{propertiesofkapparing}-\ref{npforms}. The main
tools are the syzygy bundles on $\oM_S$, introduced in \ref{syzygy},
and strong vanishing theorems they satisfy, proved in
\ref{fundamentalbundles}.

We are grateful to Alexander Polishchuk for  help with Koszul algebras,
to Gavril Farkas for help with syzygies, and
to S\'andor Kovacz, Hal Schenk, Frank Sottile, and Bernd Sturmfels for stimulating discussions.

\section{Preliminary Results}

Throughout the paper we will make use of the following
non-standard notation:
\begin{Notation} We will say that a module (or sheaf) $M$ is
{\em extended from} modules (or sheaves) with a given property if there
is a finite filtration
$$
0 = M^r \subset M^{r-1} \dots \subset M^1 \subset M^0 = M
$$
such that each of the quotients $M^i/M^{i+1}$ has the given
property. Similarly, we say that $M$ is {\it extended from}
the collection of submodules (or sheaves) $M^i/M^{i+1}$.
\end{Notation}

For a subset $T \subset S$, let
$$\pi_{T}: \oM_{S} \rightarrow \oM_{T}$$
be the tautological fibration given  by dropping the points of $S\setminus T$
(and stabilizing).

Let $s\in S$ and let $S'=S\setminus\{s\}$.
$\pi_{S'}$ has $n-1$ tautological sections, which we indicate by the elements of $S'$.
Let
$$
\psi_s := \omega_{\pi_{S'}}(S').
$$

For any subset $T\subset S$ with $|T|,|S\setminus T|\ge2$,
let $\delta_T$ be the corresponding boundary divisor of $\oM_S$.

We refer to \cite{Keel}, \cite{Kapranov93a}, and
\cite{Faber} for background material on $\oM_S$.
We will make in particular frequent use of the
formulae in
\cite{Faber} for pullbacks of tautological line bundles under
the canonical fibrations $\pi_T$.
In this section $k$ may have arbitrary characteristic.
Throughout the paper we will often abuse notation and use the
same symbol for a sheaf and its pullback under some morphism.

\begin{Notation}
For a subset $T \subset S$, let $S^T := S \setminus T$.
\end{Notation}

\begin{Lemma}\label{residues} Formation of ${\pi_{T^t}}_*(\psi_t)$
commutes with all pullbacks, and for any vector bundle $F$ on
$\oM_{T^t}$ there are canonical identifications
$$
\begin{CD}
{\pi_{T^t}}_*(\pi_{T^t}^*(F) \otimes \psi_t) = F \otimes_k W_{T^t} \\
H^0(\oM_T,\pi_{T^t}^*(F) \otimes \psi_t) = H^0(\oM_{T^t},F) \otimes_k W_{T^t}
\end{CD}
$$
induced by taking residues along the tautological sections $T^t$.

Let $[C] \in \oM_{T}$ be a stable $T$-pointed
rational curve. Taking residues at the
points $t \in T$ gives a canonical
identification
$$
H^0(C,\omega_C(T)) = W_{T}.
$$
\end{Lemma}
\begin{proof} Let $\pi:= \pi_{T^t}$. We have the natural exact sequence
$$
0 \to \omega_{\pi} \to \omega_{\pi}(\Sigma) \to \omega_{\pi}(\Sigma)|_{\Sigma} \to 0
$$
where $\Sigma$ means the disjoint union of the sections $T^t$. Taking
residues gives a canonical identification of the right hand term with
$\cO_{\Sigma}.$ Note $\psi_t = \omega_{\pi}(\Sigma)$ by definition.
$H^1(C,\omega_C(T^t)) = 0$ for a $T^t$-pointed stable curve of genus $0$,
so $R^1\pi_*(\psi_s)$ vanishes, and $R^i\pi_*(\psi_s)$ are vector bundles
and their formation commutes with all base extensions, for all $i$, by
the semi-continuity theorem.

Applying
$\pi_*$ and using the duality identification $R^1\pi_*(\omega_{\pi}) = \cO$
we
obtain a natural exact sequence
$$
0 \to \pi_*(\psi_t) \to \boplus_{x \in T^t} \cO_{\delta_{x,t}} \to \cO_{\oM_{T^t}} \to
0
$$
where the map $\cO_{x,t} \to \oM_{T^t}$ is the identity (note $\delta_{x,t}$ is
the section corresponding to $x \in T^t$). Whence the identification
$\pi_*(\psi_t) = W_{T^t} \otimes \cO$. $R^1\pi_*(\psi_t)$ is zero. The rest follows
from this and the projection formula. \end{proof}

\begin{Definition}\label{defofkappaF} For a subset
$F \subset S$, let $\kappa_F := \pi_F^*(\kappa)$
and $L_F := \kappa \otimes \kappa_F^*.$
\end{Definition}

\begin{Lemma} \label{old2.5}
Let $F \subset T \subset S$ be subsets, with $|F| \geq 3$.
Let
$$
T = T_{|T|} \subset T_{|T|+1} \dots \subset T_{|S|} = S
$$
be a flag of subsets, and define $t_i := T_i \setminus T_{i+1}$ for
$i > |T|$.
Then
$$
L_F = \pi_{T}^*(L_F) \bigotimes_{|S| \geq i > |T|} \pi_{T_i}^*(\psi_{t_i})
$$
\end{Lemma}
\begin{proof} Note by definitions
$$
L_F \otimes \pi_{T}^*(L_F) = \kappa \otimes \kappa_T^*.
$$
This we compute by induction on $|S \setminus T|$. Its enough
to consider the case $T = S^s$, and so prove the formula
\begin{equation}\label{old2.5.1}
\pi_{S^s}^*(\kappa) \otimes \psi_s = \kappa
\end{equation}
which is  given by wedge product of forms. \end{proof}

\begin{Corollary}\label{old2.7} There is a canonical
identification
$$
H^0(\oM_S,\kappa) = H^0(\oM_{S^s},\kappa_{S^s}) \otimes H^0(\oM_S,\psi_s).
$$
A flag of subsets as in Theorem \ref{maintheorem} canonically
induces identifications
$$
H^0(\oM_S,\kappa) = \bigotimes_{n \geq i \geq 4} H^0(\oM_{S_i},\psi_{s_i}) =
\bigotimes_{n \geq i \geq 4} W_{S_{i-1}}
$$
\end{Corollary}
\begin{proof} Immediate from Lemma \ref{old2.5} and Lemma \ref{residues}.
\end{proof}

\begin{Corollary}\label{old2.8}
$L_F$ is globally generated and
big. $\kappa$ is very ample. \end{Corollary}
\begin{proof} Consider first $L_F$. By \eqref{old2.5} and
induction it's
enough to consider $\psi_i$. Global generation
of $\psi_i$ is due to Kapranov, the associated map
is his birational contraction
$$
\oM_S \rightarrow \bP^{|S|-3}.
$$
See \cite[2.9]{Kapranov93a}. Now it follows from
\eqref{old2.7} and induction, that $\kappa$ is globally generated, and the
map given by global sections
$$
\oM_S \rightarrow \bP(H^0(\oM_S,\kappa)^*)
$$
factors through the Segre embedding
$$
\oM_{S^s} \times \bP(H^0(\oM_S,\psi_s)^*)
\subset \bP(H^0(\oM_{S^s},\kappa)^*) \times \bP(H^0(\oM_S,\psi_s)^*).
$$
To prove the map
$$
\oM_{S} \rightarrow \oM_{S^s} \times \bP(H^0(\oM_S,\psi_s)^*)
$$
is a closed embedding, it is enough to prove this for the restriction
to each fibre
$C = \pi_{S^s}^{-1}([C])$ of $\pi_{S^s}$. But by Lemma~\ref{residues}
restriction gives a canonical identification
$$
H^0(\oM_S,\psi_{s}) = H^0(C,\omega_C(S^s))
$$
and one checks easily that on the stable $S^s$-pointed curve $C$,
$\omega_C(S^s)$ is very ample. \end{proof}

Next we prove a topological analog of Corollary \ref{old2.7}.

\begin{Corollary}\label{specialsections}
Given $a,b \in T$ there is a
unique section
$$
\omega(a b) \in H^0(C,\omega_C(T))
$$
which has residue $1$ at $a$, $-1$ at $b$, and
is regular everywhere else. Let $a,b \in F \subset T$
with $|F| \geq 3$ and let
$$
\pi_{F}: C \rightarrow C'
$$
be the stabilisation of $(C,F)$. Pullback induces
canonical identifications
$$
H^0(C',\omega_{C'}(F)) = H^0(C,\omega_C(F))
\subset H^0(C,\omega_C(T))
$$
under which $\omega(ab)$ is sent to $\omega(ab)$. \end{Corollary}
\begin{proof} This is immediate from Lemma~\ref{residues} and the
definition of stabilization. \end{proof}

\begin{Lemma}\label{old2.2} Given distinct $a,b \in S$ there is a
a global $1$-form
$$
\omega \in H^0(\oM_S,\Omega^1(\log \partial \oM_S))
$$
whose restriction to $C \subset \oM_S$ is
$\omega(ab)$, for all $[C] \in \oM_{S^c}$.
\end{Lemma}
\begin{proof}
We have a commutative diagram
$$
\CD
\oM_S @>{\pi_{F \cup s}}>> \oM_{F \cup s} = C' = \bP^1 \\
@V{\pi_{S^s}}VV     @V{\pi_F}VV \\
\oM_{S^s} @>{\pi_{F}}>> \oM_F = {\operatorname{pt}}
\endCD
$$
where $C'$ is as in Lemma~\ref{specialsections}. Now we take for
$\omega$ the pullback of $\omega(ab)$ from $C'$.  \end{proof}

\begin{Theorem}\label{old2.3}
Let $p_1,\dots,p_{n-1}$ be distinct points
in $\bP^1$. Over the complex numbers there
are vector space identifications
$$
\begin{CD}
H^\bullet(M_n, \bC) = H^\bullet(M_{n-1},\bC) \otimes
H^\bullet(\bP^1 \setminus \{p_1,\dots,p_{n-1}\},\bC)  \\
H^\bullet(M_n, \bC) = \bigotimes_{i=3}^{n-1} H^\bullet (\bP^1 \setminus \{p_1,\dots,p_i\}, \bC)
\end{CD}
$$
and canonical identifications

$$
\begin{CD}
H^{n-3}(M_n,\bC) = \bigotimes_{i=3}^{n-1} H^1(\bP^1 \setminus \{p_1,\dots,p_i\}, \bC) \\
H^{k}(M_n, \bC) = H^0(\Omega^k(\log \partial \oM_n)) \\
H^{n-3}(M_n,\bC) = H^0(\oM_n,\kappa)
\end{CD}
$$
\end{Theorem}
\begin{proof}
For $[C] \in M_{S^s}$ and
fixed $a \in S^s$ the differential forms
$\omega(a b)$, $b \in S^{s,a}$ give a basis
of
$$
H^0(C,\omega_C(S^s)) = H^1(C \setminus S^s,\bC).
$$
By Lemma~\ref{old2.2} these are restrictions of global log
forms, and in particular global cohomology
classes. Thus by induction and
the Leray-Hirsch theorem, \cite[5.11]{BottTu}, there is
an additive isomorphism
$$
H^\bullet(\oM_S,\bC) = H^\bullet(\oM_{S^s},\bC) \otimes
H^\bullet(\bP^1 \setminus S^s,\bC).
$$
and $H^\bullet(\oM_S,\bC)$ is generated by meromorphic
$1$-forms with log poles on the boundary.
By \cite{Deligne}, such forms are never
exact. Note the last formula is just the
$k = n-3$ case of the preceding formula.
The isomorphism given by the Leray-Hirsch theorem depends
on choosing global lifts for the $\omega(a b)$. However in
the top degree it is independent of choices.
\end{proof}

\section{Filtrations}

For any globally generated line bundle $L$
on a projective variety $X$, we define a vector bundle
$V_L$ by the exact sequence
$$0\to V_L\to H^0(X,L)\otimes\cO_X\to L\to0.$$

\begin{Lemma}\label{firstVpsi} There is a natural exact sequence
$$
0 \to \pi_{S^t}^*(V_{\psi_s}) \to V_{\psi_s} \to
\cO(-\delta_{s,t}) \to 0.
$$
\end{Lemma}

\begin{Corollary}\label{filtrationofVpsi}
Choose a flag of subsets
$$
S_3 \subset S_4 \subset \dots S_n = S
$$
as in \eqref{maintheorem}.
$V_{\psi_s}$ is extended from
the line bundles
$$
\pi_{S_i}^*(\cO(-\delta_{s,s_i}))
$$
for $|S| \geq i \geq 4$.

$\wedge^q V_{\psi_s}$ is extended from line bundles of
form $\cO(-E)$ with $E$ a sum of distinct boundary
divisors with $a \equiv b \not \equiv s$, where $S_3 = \{a,b,s\}$
and we think of boundary divisors as a partition (or equivalence
relation) on $S$.
\end{Corollary}

\begin{proof}
Immediate from the Lemma, and the formula for pulling back boundary divisors under
maps $\pi_T$. \end{proof}

\begin{proof}[Proof of Lemma~\ref{firstVpsi}]
We have
$$
\psi_s = \pi_{S^t}^*(\psi_s)\otimes\cO(\delta_{s,t})
$$
and $\psi_s|_{\delta_{s,t}}$ is canonically trivial
(by taking residues). This induces a commutative diagram
with exact rows and columns
$$
\begin{CD}
  @. 0 @. 0 @. 0 \\
@. @VVV @VVV @VVV @. \\
0 @>>> \pi_{S^t}^*(V_{\psi_s}) @>>> V_{\psi_s} @>>> \cO(-\delta_{s,t}) @. \\
@. @VVV @VVV @VVV @. \\
0 @>>> \pi_{S^t}^*(H^0(\psi_s))\otimes_{k} \cO_{\oM_S}
@>>> H^0(\psi_s)\otimes_{k} \cO_{\oM_S} @>>> H^0(\cO_{\delta_{s,t}})
\otimes_k \cO_{\oM_S}
@>>> 0 \\
@. @VVV @VVV @VVV @. \\
0 @>>> \pi_{S^t}^*(\psi_s) @>>> \psi_s @>>> \cO_{\delta_{s,t}} @>>> 0 \\
@. @VVV @VVV @VVV @. \\
 @. 0 @. 0 @. 0 \\
\end{CD}
$$
Here the second row is obtained from the third by taking global
sections -- which gives a short exact sequence of vector spaces
as the $H^1$ term vanishes -- and then tensoring with
the structure sheaf $\cO_{\oM_S}$. Now the first row is given
by taking kernels of the vertical maps.

Now the result follows from the snake lemma. \end{proof}

\begin{Lemma}\label{boundarylemma} Let
$a,b,s \in S$ be distinct. Let $E$ be the (reduced) union
of all boundary divisors of $\oM_S$ with
$a \equiv b \not \equiv s$. Then there is member of
$|\kappa|$ which is an effective combination of irreducible
components of $B - E$.
\end{Lemma}

\begin{proof} We induct on $|S|$. For $|S| =4$ we can take
$\delta_{a,s}$. By induction it is enough to find a member
of $|\psi_t|$, $t \in S^{a,b,s}$, supported on components
of $B - E$. This is clear from Lemma~\ref{memberofpsi}. \end{proof}

\begin{Lemma}\label{memberofpsi}
$$
\sum_{a,s \in F,t \in F^c} \delta_{F}
$$
is linearly equivalent to $\psi_t$.
\end{Lemma}
\begin{proof} The sum above is
$\Psi_t^*(\Psi_t(\delta_{a,s}))$ where
$$
\Psi_t: \oM_S \rightarrow \bP^{|S| -3}
$$
is the Kapranov model. $\Psi_t(\delta_{a,s})$ is
a hyperplane, and so the above pullback represents
$\psi_t$. \end{proof}

\begin{Lemma}\label{mainfiltration}
$\wedge^q V_{\psi_s} \otimes \kappa$ is extended from
globally generated line bundles. It is also extended from
line bundles of form $\cO(E)$ where $E$ is a divisor
linearly equivalent to a $\bQ$-divisor of form
$K_{\oM_S} + \Delta + A$, for $\Delta$ an effective
combination of boundary divisors with coefficients
strictly less than one, and $A$ an ample divisor.
\end{Lemma}
\begin{proof} By Corollary~\ref{filtrationofVpsi},
$\wedge^q V_{\psi_s} \otimes \kappa$
is extended from line bundles associated to divisors
$$
\kappa + \sum_{t=1}^{t=r} \pi_{S_{i_t}}^*(\delta_{s,s_{i_t}})
$$
for some sequence of integers
$$
n \geq r > r-1 \dots > r_1 \geq 4.
$$
Global generation of such a divisor follows by induction, the formulae
$$
\kappa_S = \kappa_{S^t} + \psi_t
$$
$$
\psi_t = \pi_{S^s}^*(\psi_t) + \delta_{s,t}
$$
and the global generation of $\kappa$ and $\psi_i$.

The second statement is immediate from
Lemma~\ref{filtrationofVpsi} and Lemma~\ref{boundarylemma}.
\end{proof}
\begin{Corollary}\label{mainvanishing}
$H^i(\wedge^q V_{\psi_s} \otimes \kappa \otimes M) = 0$ for
any $i > 0$ and any nef line bundle~$M$.
\end{Corollary}
\begin{proof} This is immediate from the preceding
Lemma and the Kawamata-Viehweg vanishing theorem. \end{proof}

\begin{Corollary}\label{npforVeronese}
Let $(C,S^s)$ be an $S^s$-pointed stable curve and $\psi := \omega_C(S^s)$.
Then $H^i(C,\wedge^q (V_{\psi}) \otimes \psi) = 0$ for $i > 0$,
and $H^0(\wedge^q(V_{\psi})) = 0$ for all $q > 0$.
The log canonical embedding $C \subset \bP(H^0(C,\omega_C(S^s))^*)$
satisfies Green--Lazarsfeld's condition $N_p$ for all $p \geq 0$.
\end{Corollary}
\begin{proof}
$C = \pi_{S^s}^{-1}[C] \subset \oM_S$, $\psi = \psi_s|_C = \kappa|_C$, and
$V_{\psi} = {V_{\psi_s}}|_C.$ $\wedge^q V_{\psi_s}$ is extended from line
bundles satisfying the Kawamata-Viehweg vanishing theorem by
Lemma~\ref{mainfiltration}, and restriction.
This gives the vanishing result, which implies $N_p$ for all $p$, see
\cite{GreenLazarsfeld}. By
Lemma~\ref{filtrationofVpsi} and restriction, $\wedge^q (V_{\psi})$ is
extended from line bundles with no global sections, and thus
has itself no global sections.
\end{proof}

\section{The syzygy bundles}\label{syzygy}

We begin by giving a canonical resolution of the
structure sheaf
of $\oM_S$ by natural vector bundles over $\oM_{S^s} \times \bP^{n-3}$,
using a special case of the Beilinson spectral sequences as in
\cite{GrusonLazarsfeldPeskine}.

The diagonal embedding
$$
\bP^r \subset \bP^r \times \bP^r
$$
is the zero locus of a regular tautological section
of $V_{\cO(1)} \boxtimes \cO(1)$.
For any morphism $f: X \rightarrow \bP^r$ the section pulls back to a
regular section of $V_{\cL} \boxtimes \cO(1)$, $\cL := f^*(\cO(1))$,
on $X \times \bP^r$ with
zero locus the graph $X = \Gamma_f \subset X \times \bP^r$. The Koszul
complex then gives an exact sequence
\begin{equation}\label{Koszulcomplex}
0 \to \wedge^{r} V_{\cL} \boxtimes \cO(-r) \to \dots \to V_{\cL} \boxtimes
\cO(-1) \to \cO_{X \times \bP^{r}} \to \cO_X \to 0.
\end{equation}

We have a closed
embedding $\Phi: \oM_S \subset \oM_{S'} \times \bP^{n-3}$ given
by $\pi := \pi_{S^s}$ and the linear series $|\psi_s|$, see the
proof of Corollary~\ref{old2.8}.

\begin{Lemma}\label{defofMq}
$\cM^q := R^1\pi_*(\wedge^{q+1} V_{\psi_s})$ is a vector bundle on $\oM_{S^s}$
for $q \geq 0$. $\cM^0 = 0$.
One has an exact sequence of sheaves
on $\oM_{S'}\times\bP^{N-3}$:
\begin{equation}\label{sheafcomplex}
0\to\cM^{N-4}\boxtimes\cO(3-N)\to\ldots\to\cM^1\boxtimes\cO(-2)\to\cO_{\oM_{S'}\times\bP^{N-3}}
\to\Phi_*\cO_{\oM_S}\to0
\end{equation}
The fibre of $\cM^q$ at the point $[C] \in \oM_{S^s}$ is canonically identified
with the $q$-th syzygy for the $S^s$-pointed stable curve $(C,S^s) \subset \bP^{n-3}$, i.e.
$$
\cM^q|_{[C]} = \Tor_q^A(k,B)
$$
where $A = R(\bP^{n-3},\cO(1))$ is the homogeneous coordinate ring of $\bP^{n-3}$ and
$B= R(C,\omega_C(S^s))$ is the homogeneous coordinate ring of $C \subset \bP^{n-3}$,
for $q > 0$.
\end{Lemma}

\begin{proof}
We apply the above construction in the case of $X = \oM_S \rightarrow
\bP^{n-3}$, and
push forward the sequence \eqref{Koszulcomplex} along
$$
p=\pi \times id: \oM_S \times \bP^{n-3}
$$
where
the fibres of $p$ (or equivalently $\pi$) are $S^s$-pointed stable curves. So
the formation of $R^i\pi_*(\wedge^{q+1} V_{\psi_s})$ commute with all
base extensions, vanish for $i =0$ or $i =q =1$, and are vector bundles
for $i=1$, by Lemma~\ref{npforVeronese}.

Exactness of \eqref{sheafcomplex} follows by analyzing
the spectral sequence for the hyperderived pushforward:
$$
E^1_{i,j} = R^jp_*(\cF^i),\quad i \leq 1,\quad j \geq 0
$$
with $\cF^{-i} = \wedge^i(V_{\psi_s}) \boxtimes \cO(-i)$,
$i \leq 0$, $\cF^{1} = \cO_{\oM_S}$. Since \eqref{Koszulcomplex}
is exact, the spectral sequence abuts to zero. By
\ref{npforVeronese}, the sequence has only two non-zero rows:
$$
E^1_{*,1}: \dots \to R^1p_*(\wedge^3 V_{\psi_s}) \boxtimes \cO(-3) \to
R^1p_*(\wedge^2 V_{\psi_s}) \boxtimes \cO(-2) \to 0 \to 0 \dots
$$
and
$$
E^1_{*,0}: \dots 0 \to 0 \to \cO_{\oM_{S^s} \times \bP^{n-3}}
\to \Phi_*(\cO_{\oM_S}) \to 0 \to 0 \dots.
$$
The exactness of \eqref{sheafcomplex} follows easily.

If we restrict \eqref{sheafcomplex} to $[C] \times \bP^{n-3}$ it remains exact
(for example by the flatness of~$\pi$) and we obtain the exact sequence of
sheaves on $\bP^{n-3}$:
$$
0\to H^1(C,\wedge^{n-3} V_{\psi}) \otimes \cO_{\bP^{n-3}}(-(n-3))
\to \dots \to H^1(C,\wedge^2 V_{\psi}) \otimes \cO(-2) \to \cO_{\bP^{n-3}}
\to \cO_C \to 0
$$
where $\psi := {\psi_s}|_C$ as in Lemma~\ref{npforVeronese}.
Tensor with $\boplus_{n\geq0} \cO(n)$ and take global sections, a simple
spectral sequence analysis shows the resulting sequence of $R(\bP^{n-3},\cO(1))$
modules is exact. This yields a resolution of the homogeneous coordinate ring
of $R(C,\psi)$. Note each of the maps is given by linear forms, and thus when
we tensor with $k$, the maps in the resulting complex are all zero, and thus
the terms of the sequence are identified with the relevant
$\operatorname{Tors}$. \end{proof}

Let $T = S^s$. Observe the vector bundle, $\cM^i$, on $\oM_{T}$
is intrinsic to $\oM_{T}$ (it does not depend on $S$).
When there is the possibility of confusion we
refer to it as $\cM^i_{T}$.

\begin{Lemma}\label{pushpull} Let $s,t \in S$.
Let $E$ be a vector bundle on $\oM_{S^t}$.
There are canonical identifications
$$
R^p{ \pi_{S^s}}_*(\pi_{S^t}^*(E)) =
\pi_{S^{s,t}}^*(R^p{\pi_{S^{s,t}}}_*(E)).
$$
for all $p \geq 0$.
\end{Lemma}
\begin{proof}
Consider first the commutative pull back diagram:
$$
\begin{CD}
\oM_{S^s} \times_{\oM_{S^{s,t}}} \oM_{S^t} @>{\pi_2}>> \oM_{S^t} \\
@V{\pi_1}VV @V{\pi_{S^{s,t}}}VV \\
\oM_{S^s} @>{\pi_{S^{t,s}}}>> \oM_{S^{s,t}}
\end{CD}
$$
$$
R^p{ \pi_1}_*(\pi_2^*(E)) =
\pi_{S^{t,s}}^*(R^p{\pi_{S^{s,t}}}_*(E))
$$
since $\pi_{S^{s,t}}$ is flat,
see \cite[9.3]{Hartshorne}.
Now by the Leray spectral sequence and the projection formula it
is enough to show that
$$
R^qf_*(\cO) =0
$$
for $q > 0$ where
$$
f: \oM_{S} \rightarrow \oM_{S^s} \times_{\oM_{S^{s,t}}} \oM_{S^t}
$$
is the natural map. This holds as the map is birational and
the image has rational singularities, see
e.g. \cite[pg. 548]{Keel}. \end{proof}

\begin{Theorem}\label{b2.1} For $t \in A$ there is a vector bundle
$Q$ on $\oM_A$ and exact sequences
\begin{equation}\label{b2.1.1}
0 \to \cM^p_{A^t} \to \cM^p_{A} \to Q \to 0
\end{equation}
\begin{equation}\label{b.2.1.1}
0 \to \wedge^p(V_{\psi_t}) \to Q \to \cM^{p-1}_{A^t} \to 0
\end{equation}
\end{Theorem}
\begin{proof} Let $A = S^s$.
The exact sequence \eqref{firstVpsi} yields
$$
0 \to
\pi_{S^t}^*(\wedge^{p+1}(V_{\psi_s})) \to \wedge^{p+1}(V_{\psi_s}) \to
\pi_{S^t}^*(\wedge^p(V_{\psi_s}) \otimes \cO(-\delta_{s,t})
\to 0
$$

${\pi_{S^s}}_*(\cO(-\delta_{s,t})) =0$ since $\delta_{s,t}$
is a section of $\pi_{S^s}$. So by Lemma \ref{pushpull} we
obtain the exact sequence
$$
0 \to \pi_{S^{s,t}}^*(\cM^p) \to \cM_{S^s}^p \to Q \to 0
$$
where
$$
Q := R^1{\pi_{S^s}}_*(\pi_{S^t}^*(\wedge^p(V_{\psi_s})) \otimes \cO(-\delta_{s,t})).
$$
We study $Q$ beginning with the exact sequence:
$$
0 \to \pi_{S^t}^*(\wedge^p(V_{\psi_s})) \otimes \cO(-\delta_{s,t}) \to \pi_{S^t}^*(\wedge^p(V_{\psi_s}))
\to \pi_{S^t}^*(\wedge^p(V_{\psi_s}))|_{\delta_{s,t}} \to 0.
$$
Since $\delta_{s,t}$ is a section of $\pi_{S^s}$, $R^1{\pi_{S^s}}_*$ vanishes
on the right term. ${\pi_{S^s}}_*$ applied to the middle term gives
$$
{\pi_{S^s}}_*(\pi_{S^t}^*(\wedge^p(V_{\psi_s}))) =
\pi_{S^{t,s}}^*({\pi_{S^t}}_*(\wedge^q(V_{\psi_s}))) = 0
$$
by Lemma \ref{pushpull} and Lemma \ref{npforVeronese}.
So we obtain the exact sequence:
$$
0 \to \wedge^p(V_{\psi_t}) \to Q \to \pi_{S^{s,t}}^*(\cM^{p-1}_{S^{s,t}})
\to 0.
$$
\end{proof}
We have the immediate:

\begin{Corollary}\label{filtration} A flag of
subsets
$$
S_3 \subset S_4 \subset \dots \subset S_n = S
$$
as in \eqref{maintheorem} extends
$\cM^p_S$ from vector bundles of form
$\pi_{S_i}^*(\wedge^q V_{\psi_{s_i}})$ for $p \geq q$.
\end{Corollary}

\begin{Corollary}\label{c1Mp}
$$
c_1(\cM^p_S) = { {|S|-4} \choose {p-1} }  \kappa
$$
\end{Corollary}
\begin{proof} An easy induction argument using Theorem \ref{b2.1} and
Lemma \ref{old2.7}.
\end{proof}

\begin{Remark}
As $\cM^k_T$ is the bundle with fibre at $[C]$ the $k$-th syzygy module, it
is naturally a subbundle of the trivial bundle
$\Sym_2(\bV) \otimes \bV^{\otimes(k-1)}$
where $\bV$ is the trivial bundle $\bV := \pi_*(V_{\psi_s})$, e.g.
$\cM^1 \subset \Sym_2(\bV) \otimes \cO$ has fibre at $[C]$ the
space of conics vanishing on $C \subset \bP^{|T| - 2}$,
$\cM^2 \subset \cM^1 \otimes \bV$ has fibre the linear relations
among these quadrics, etc.
As a sub bundle of a trivial bundle,
$\cM^k$ induces a map of $\oM_{T}$ to a Grassmannian. By Corollary \ref{c1Mp}
the first Chern class of $\cM^k$ is ample, thus the map is finite, and one
can check (with a bit of work) that it is in each case a closed embedding,
(at least) two of which have been previously studied:
As noted $\cM^1 \subset \Sym_2(\bV)$ gives fibrewise the space
of conics vanishing on $C \subset \bP^{|T| -2}$.
By Lemma \ref{npforVeronese} $C$ is cut out by such quadrics. The
induced closed embedding is Kapranov's realization
of $\oM_{S^s}$ as the closure in the Hilbert scheme of the locus of rational
normal curves in $\bP^{|T|-2}$ through $|T|$ fixed points. See \cite{Kapranov93a}.
Let $T := S^s$, $\pi := \pi_{T}$.
$$
\cM^{|T|-3} = R^1\pi_*(\wedge^{|T|-2} V_{\psi_s}) = R^1\pi_*(\omega_\pi(\Sigma))
$$
where $\Sigma$ is the union of the $|T|$ tautological sections. It follows
then by the deformation theory of pairs that
$\cM^{|T|-3}$ is the log tangent bundle $T_{\oM_{T}}(-\log B)$,
and one can check that the map to the Grassmannian is given by the space of global log $1$-forms.
The induced map to projective space (composing with the Pl\"ucker
embedding of the Grassmannian) is the log canonical embedding of $M_{T}$,
and is an instance of a general construction that holds for any complement
to a hyperplane arrangement, see \cite{HackingKeelTevelev}.
\end{Remark}

\section{Koszulness and $p$-linearity}

In this section $A = \boplus_{n \geq 0} A_n$ is a commutative
graded algebra over a field $A_0 = k$ and $M$ is a graded $A$-module.
All our modules are concentrated in nonnegative degrees.

\begin{Definition}
$M$
is called {\em $p$-linear} if
$\Tor_i^A(M,k)$ is concentrated in
degrees $i,i+1,\ldots,i+p$.
$A$ is called {\em Koszul} if the trivial $A$-module $A_0=k$ is $0$-linear.
\end{Definition}

\begin{Example}\label{Poly}
A polynomial algebra $S$ over $k$ (graded by total degree) is Koszul:
the minimal resolution of $k$ over $S$ is given by the standard Koszul complex.
\end{Example}

Among the many nice properties of Koszul algebras we note the
result of \cite[2.3.3]{BeilinsonGinszburgSoergel}:
\begin{Theorem}\label{Koszulisquadratic}
Let $A$ be a Koszul algebra. Then the
natural map
$$
\Sym_{A_0}(A_1) \rightarrow A
$$
is surjective, and its kernel is generated by degree two
elements.
\end{Theorem}

\begin{Lemma}\label{torsequence}
Consider
a short exact sequence of graded $A$-modules:
\begin{equation}\label{shortseq}
0\to M'\to M\to M''\to0.
\end{equation}
If $M'$ and $M''$ are $k$-linear then $M$ is $k$-linear.
If $M$ is $k$-linear and $M''$ is $(k-1)$-linear then $M'$ is $k$-linear.
If $M$ is $k$-linear and $M'$ is $(k+1)$-linear then $M''$ is $k$-linear.
\end{Lemma}

\begin{proof}
All statements immediately follow from the long exact sequence for $\Tor$.
\end{proof}

\begin{Corollary}\label{filtr}
If $M$ is extended from $k$-linear modules
then $M$ is $k$-linear.
\end{Corollary}

\begin{proof}
Induction on the depth of the filtration using Lemma~\ref{torsequence}.
\end{proof}

\begin{Lemma}\label{nilpotent}
If $A$ is Koszul and $M$ is concentrated in degrees at most $p$
then $M$ is $p$-linear.
\end{Lemma}

\begin{proof}
$M$ has the natural filtration
$$
M = M_{\geq 0} \supset M_{\geq 1} \supset \dots \supset M_{\geq p+1} = 0
$$
where $M_{\geq i} := \boplus_{n \geq i} M_n$. Thus $M$ is
extended from the modules $M_{\geq i}/M_{\geq i+1}$.
$M_{\geq i}/M_{\geq i+1}$ is a $k$-module (i.e.
annihilated by $A_n$, $n > 0$) concentrated in degree $i$,
and thus $i$-linear by the definition of a Koszul
algebra. Thus $M$ is $p$-linear by Corollary \ref{filtr}.
\end{proof}

\begin{Lemma}\label{truncation}
If $A$ is Koszul then
$M$ is $k$-linear iff $M':=\boplus_{n\ge k}M_n$ is $k$-linear.
\end{Lemma}

\begin{proof}
Consider the short exact sequence \eqref{shortseq}, where
$M''=\oplus_{n<k}M_n$ has
a natural structure of a nilpotent $A$-module.
$M''$ is $(k-1)$-linear by Lemma~\ref{nilpotent}.
Therefore $M$ is $k$-linear by Lemma~\ref{torsequence}.
\end{proof}

\begin{Lemma}\label{funny}
Consider the complex of graded $A$-modules
\begin{equation}\label{complex1}
0 \to M^r\to
\ldots \to M^3\to M^2\to M^1\to M^0 \to 0.
\end{equation}
Assume that $M^k$ is $p+k$-linear, for $k \geq 1$,
and that the cohomology module $H^k$ is $p+k+1$-linear
for $k \geq 0$. Then $M^0$ is $p+1$-linear.
\end{Lemma}

\begin{Remark}
We will use this Lemma in the situation when
$A$ is Koszul and $H^k$ is concentrated in degrees at most $p+k+1$, and
thus $p+k+1$-linear by Lemma~\ref{nilpotent}.
\end{Remark}

\begin{proof}
Let $d^k:\,M^k\to M^{k-1}$ be the differential
with kernel $\Ker^k$ and image $\Im^k$.

We argue by induction, the case $r=0$ being obvious.
Consider the complex
$$
0 \to M^r\to \ldots \to M^3\to M^2\to \Ker^1 \to 0
$$
that obviously has the same cohomology as the original complex.
It follows by induction that $\Ker^1$ is $(p+2)$-linear.

By Lemma~\ref{torsequence} and the sequence
$$0\to\Ker^1\to M^1\to\Im^1\to0$$
it follows that $\Im^1$ is $(p+1)$-linear.

By Lemma~\ref{torsequence} and the sequence
$$0\to\Im^1\to M^0\to H^1\to0$$
it follows that $M^0$ is $(p+1)$-linear.
\end{proof}

\begin{Definition}[Segre Product]
For graded $k$-modules $N$ and $M$ let
$$
N \segp M := \boplus_{n \geq 0} N_n \otimes_k M_n
$$
\end{Definition}

\begin{Proposition}\label{Segretricks}
Let $A$ and $B$ be Koszul algebras.
Then $A \segp B$ is a Koszul algebra.
Moreover, if $M$ is $p$-linear
over $A$, and $N$ is $p$-linear over $B$, then
$M \segp N$ is $p$-linear over $A \segp B$.
\end{Proposition}

\begin{proof}
For $p=0$ this is the content of \cite[6.5]{CHTV}. In particular
(taking $M = N = k$), $A \segp B$ is Koszul.

For a graded vector space $V$ let
$$
V\langle p\rangle := \boplus_{n \geq 0} V_{n+p}.
$$
Observe $(M \segp N)\langle p\rangle = M\langle p\rangle \segp N\langle p\rangle$, and that
by Lemma~\ref{truncation} $M$ is $p$-linear iff $M\langle p\rangle$ is
$0$-linear. Thus the $p$-linear case is reduced to the
zero linear case. \end{proof}

The following is a slight generalisation of
\cite[Theorem 1.2]{Polishchuk}.

\begin{Lemma}\label{genPolishchuk} Let
$A\to B$
be a homomorphism of graded rings, with $A_0 = B_0 = k$. Assume $B$ is
$1$-linear over~$A$. Let
$M$ be a graded $B$-module. If $M$ is $p$-linear over~$A$, then $M$ is
$p$-linear over~$B$.
\end{Lemma}
\begin{proof} The case of $M=k$ is \cite[1.2]{Polishchuk}. The same
proof works for any $M$. \end{proof}

\begin{Definition}
For any sheaves $\cL$, $\cM$ on an algebraic variety $X$,
we define
$$\cG^X(\cM,\cL)=\boplus_{n\ge0}\cG^X_n(\cM,\cL),\quad \hbox{\rm where}
\quad \cG^X_n(\cM,\cL)=H^0(X,\cM\otimes\cL^{\otimes n}).$$
We drop $X$ from notation if it's clear from the context.
We let $\cG(\cL):=\cG(\cO_X,\cL)$.
Notice that $\cG(\cL)$ is a graded $k$-algebra and $\cG(\cM,\cL)$ is a graded $\cG(\cL)$-module.
We call $\cL$ a {\em Koszul sheaf} if $\cG(\cL)$ is a Koszul algebra.
We call $\cM$ {\em $p$-linear over} $\cL$ if $\cG(\cM,\cL)$ is $p$-linear over $\cG(\cL)$.
\end{Definition}

\begin{Lemma}\label{Np} Let $\cL$ on $X$ be a very ample line bundle
and assume the coordinate ring $B := \cG^X(\cL)$ is Koszul.
Assume the
embedding $f:\,Y\hookrightarrow X$ is non-degenerate, i.e.
$$
H^0(X,\cL) \to H^0(Y,\cL|_Y)
$$
is an isomorphism.
Then $f_*(\cO_Y)$ is $1$-linear over $\cL$ iff the minimal resolution of
$R := \cG^Y(\cL|_Y)$ over $B$ is of form
$$
\dots \to B[-k]^{a_k} \to \dots \to B[-2]^{a_2} \to B \to R \to 0
$$
\end{Lemma}
\begin{proof}
Clearly if we have such a resolution then $R$ is $1$-linear. So assume
$R$ is $1$-linear over $B$. By Lemma \ref{genPolishchuk} $R$ is Koszul, thus
generated by degree $1$ elements by Theorem \ref{Koszulisquadratic}.
Thus $B \to R$ is surjective and the kernel, $K$, is generated
by elements of degree at least $2$. Consider now a minimal free
resolution of $K$:
$$
\dots \to  F_k \to \dots \to F_2 \to K \to 0.
$$
Clearly minimal generators of $F_k$ have degree at least $k$. But $K$ is
$2$-linear by Lemma \ref{torsequence}, thus minimal generators of $F_k$
are of degree exactly $k$, see e.g. the second paragraph of
\cite{ParanjapeRamanan}. \end{proof}

\begin{Lemma}\label{Segre}
Let $\cL_i$, $\cM_i$ be sheaves on $X_i$, $i=1,2$.
Assume that $\cL_i$ is Koszul and $\cM_i$ is $p$-linear over $\cL_i$, $i=1,2$.
Then $\cL_1\boxtimes\cL_2$ is a Koszul sheaf on $X_1\times X_2$
and $\cM_1\boxtimes\cM_2$ is $p$-linear over  $\cL_1\boxtimes\cL_2$.
\end{Lemma}

\begin{proof}
Immediate from Lemma~\ref{Segretricks}.
\end{proof}

\begin{Lemma}\label{sheafpoli}
Let $f:\,Y\hookrightarrow X$ be a closed embedding.
Let $\cL$ be the Koszul sheaf on $X$.
Assume that $f_*(\cO_Y)$ is $1$-linear over $\cL$.
Then $\cL|_Y$ is Koszul.
Moreover, if $\cM$ is the sheaf on $Y$ and $f_*(\cM)$ is $p$-linear
over $\cL$ then $\cM$ is $p$-linear over $\cL|_Y$.
\end{Lemma}

\begin{proof}
Immediately follows from Lemma~\ref{genPolishchuk}.
\end{proof}

\section{Fundamental vector bundles}\label{fundamentalbundles}
\begin{Definition}
By a {fundamental} vector
bundle on $\bP^m$ we mean a bundle of form $\Lambda^q V$
for $q \geq 0$ where $V$ is the universal rank
$m$ subbundle.

By a {\it fundamenal} vector bundle on $\oM_S$ we mean the
pullback of a fundamental vector bundle on projective space
under a composition $\Psi_t \circ \pi_T$
for $t \in T \subset S$, where
$\Psi_t: \oM_T \to \bP^{|T| -3}$ is Kapranov's map.
\end{Definition}

For any Young diagram $\lambda$ with at most $N-3$ rows, let $S_\lambda$ be the corresponding
Schur functor and let $S_\lambda(V)$ be the corresponding
vector bundle on $\bP^{N-3}$.
For example, vector bundles $\Lambda^kV$ correspond to $k$-box diagrams with $k$ rows
(of one box each).

\begin{Lemma}\label{decompose}
Let $D$ be the tensor product of at most $p$ fundamental vector bundles on $\bP^{m}$.
Then $D$ is the direct sum of vector bundles of the form $S_\lambda(V)$,
where all rows of $\lambda$ have at most $p$ boxes.
\end{Lemma}

\begin{proof}
Follows from the Littlewood--Richardson rule.
\end{proof}

\begin{Proposition}\label{main}
If $\lambda$ has at most $p$ boxes in each row then
$$H^j(S_\lambda(i))=0$$
for $j>0$, $i\ge p-j$.
\end{Proposition}

\begin{proof}
By the Borel--Bott--Weyl theorem, see \cite[8.1]{Tevelev} the following holds:
Suppose
$$
\lambda = (\lambda_1,\lambda_2,\dots,\lambda_m)
$$
with $\lambda_1 \geq \dots \geq \lambda_m \geq 0$. If
for some $i$, $r = \lambda_i - i$, then all the cohomologies
vanish (the so called singular case of the theorem). Otherwise
there is a unique $j$ such that
$$
H^j(S_{\lambda}(r)) \neq 0
$$
and it is described as follows: If $r > \lambda_1 -1$ then $j =0$.
If $\lambda_m - m > r$ then $j = m$. Otherwise $j$ is the
unique index $i$ so that
$$
\lambda_i - i > r > \lambda_{i+1} - (i+1).
$$

In particular, if $r \geq \lambda_i -i$ then $H^i(S_{\lambda}(r)) = 0$.
Since $p \geq \lambda_i$, the result follows. \end{proof}

\begin{Lemma}\label{Pvanish}
Let $D$ be the tensor product of at most $p$ fundamental vector bundles on $\bP^{m}$.
Then $D(n)$ has no higher cohomology for $n\ge p$.
\end{Lemma}

\begin{proof}
Follows from Lemma~\ref{decompose} and Proposition~\ref{main}.
\end{proof}

\begin{Lemma}\label{cohomology}
If $\tilde D$ is a product of at most $p+1$ fundamental vector bundles on $\oM_{S'}$
then $H^i(\tilde D\otimes(\kappa')^{\otimes n})=0$ for any $i>0$ and $n\ge p+1$.
\end{Lemma}

\begin{proof} Immediate from the Kawamata-Viehweg vanishing theorem
and Lemma~\ref{mainfiltration}. \end{proof}

\begin{Corollary}\label{cohomologyvanish}
If $D$ is a product of at most $p$ fundamental vector bundles on $\oM_{S}$
then $H^i(\cM^k\otimes D\otimes\kappa^{\otimes n})=0$ for any $i>0$, $k>0$, and $n\ge p+1$.
\end{Corollary}

\begin{proof}
Follows from Lemma~\ref{cohomology} and Proposition~\ref{filtration}.
\end{proof}

\begin{Proposition}\label{LinearityD''} Let
$D$ be the tensor product of at most $p$ fundamental
vector bundles on $\bP^m$, $p \geq 0$ (when $p=0$,
$D = \cO$). Then $D$ is $p$-linear
over $\cO_{\bP^m}(1)$.
\end{Proposition}
\begin{proof}
By Lemma~\ref{truncation}, it suffices to prove that
$$M':=\boplus_{n\ge p}H^0(\bP^{N-3},D(n))$$
is $p$-linear over $\cO(1)$.
$D$ is the tensor product of at most $p$ fundamental vector bundles of the form $\Lambda^{n_i}V$.
We argue by induction on $\sum n_i$.

If all $n_i=0$ then $p=0$ and $D = \cO$, which is obviously $0$-linear.
Otherwise suppose we have a tensor factor $\Lambda^q V$ with $q > 0$.
Wedging the defining sequence for $V$ gives an exact sequence
$$
0 \to \Lambda^q V \to E \to (\Lambda^{q-1}V) \otimes \cO(1) \to 0
$$
with $E$ a trivial bundle. Thus by Lemma~\ref{Pvanish}
the module $M'$
sits in an exact sequence
$$0\to M'\to M\to M''(1)\to0,$$
where $M$ and $M''$ are products of trivial vector bundles with
products of at most $p$ fundamental bundles with smaller $\sum n_i$.
Now the result follows from Lemma~\ref{torsequence} and induction.
\end{proof}

\section{Koszulness of $\kappa$} \label{kappa}

We fix a filtration
$$
S_4 \subset S_5 \dots \subset S_N = S
$$
as in Theorem \ref{maintheorem}. Let $S' = S_{N-1}$.

Let $\pi:\,\oM_S\to\oM_{S'}$ be the projection map,
let $\Psi:\,\oM_S\to\bP^{N-3}$ be Kapranov's map.
Let $\kappa'$ be the log canonical line bundle on $\oM_{S'}$.
We have a closed embedding
$$
\Phi: \oM_S \subset \bP:= \bP^1 \times \bP^2 \dots \times \bP^{n-3}.
$$
Define
$$
\cL := \cO_{\bP^1}(1) \boxtimes \cO_{\bP^2}(2) \dots \boxtimes \cO_{\bP^{n-3}}(1).
$$
\begin{Theorem}
$\kappa$ is Koszul. $\cO_{\oM_S}$ is $1$-linear over $\cL$.
\end{Theorem}
\begin{Remark} Thus $\oM_S \subset \bP$ satisfies the analog of the
Green-Lazarsfeld condition  $N_p$
for all $p$, see Lemma ~\ref{Np}.\end{Remark}
\begin{proof}
We argue by induction on $N$, the case of $\oM_{0,4}=\bP^1$ being obvious.
$\cL$ is a Koszul sheaf by Example~\ref{Poly} and
Lemma~\ref{Segretricks}. So by Lemma~\ref{genPolishchuk}
its enough to show $\cO_{\oM_S}$ is $1$-linear over $\cL$.
For this we prove simultaneously by induction

\begin{Theorem}\label{inductiveversion}
Let $D$ be the tensor product of at most $p$ fundamental vector bundles
(if $p=0$ then let $D=\cO$) on $\oM_S$.
Then $\Phi_*(D)$ is $(p+1)$-linear over $\cL$.
\end{Theorem}

It's clear that $D$ can be written uniquely as
$$D=\pi^*(D')\otimes\Psi^*(D''),$$
where $D'$ is a tensor product of at most $p$ fundamental vector bundles
on $\oM_{S'}$ and $D''$ is a tensor product of at most $p$ fundamental vector bundles
on $\bP^{N-3}$ of the form $\Lambda^kV$, where $V$ is the tautological quotient bundle.

For any sheaf $\cM$ on $\bP$, let
$$\cD(\cM)=\boplus_{n\ge p+1}\cG_n(\cM,\cL).$$
The exact sequence of sheaves \eqref{sheafcomplex} induces the complex of $\cG(\cL)$-modules

\begin{equation}\label{BigComplex}
\ldots\to\cD([D'\otimes\cM^1]\boxtimes D''(-2))
\to\cD(D'\boxtimes D'')\to\cD(\Phi_*(D)) \to 0
\end{equation}

Let
$\bP' := \bP^1 \times \bP^2 \dots \times \bP^{n-4}$
and let $\cL'$ be the analog of $\cL$ on $\bP'$.

\begin{Lemma}\label{linearityD'}
$D'\otimes\cM^q$ is $(p+2)$-linear over $\cL'$ for all $q \geq 1$.
$D'$ is $p+1$-linear over
$\cL'$.
\end{Lemma}

\begin{proof}
We treat the case of $D' \otimes \cM^q$, the argument for $D'$ is entirely
analogous.
By Lemma~\ref{truncation}, it suffices to prove that
$M:=\boplus_{n\ge p+2}\cG_n(D'\otimes\cM^q,\kappa')$ is $(p+2)$-linear
over $\cG(\cL')$. By Proposition~\ref{filtration} and
Corollary~\ref{cohomologyvanish}, $M$ is extended from
modules of the form $\boplus_{n\ge p+2}\cG_n(\tilde D,\kappa')$,
where $\tilde D$ is a product of at most $p+1$ fundamental vector bundles.
So the statement follows from Corollary~\ref{filtr}
and the inductive assumption about $\oM_{S'}$.
\end{proof}

\begin{Proposition}\label{BBW}
The $i$-th cohomology module of \eqref{BigComplex} lives in degrees $\le p+i$.
\end{Proposition}
\begin{proof}
Tensoring the exact sequence \eqref{sheafcomplex} with
$[D'\boxtimes D'']\otimes\boplus_{n\ge p+1}\cL^{\otimes n}$
 induces an exact sequence of graded sheaves
\begin{equation}\label{GradedComplex}
\ldots\to \cF^2\to\cF^1\to\cF^0 \to 0,
\end{equation}
where (for $k \geq 2$)
$\cF^k =[D'\otimes \cM^{k-1}\boxtimes D''(-k)]\otimes \boplus_{n\ge p+1}\cL^{\otimes n}$.
Notice that the complex \eqref{BigComplex} is the complex of global sections of \eqref{GradedComplex}.
Since \eqref{GradedComplex} is exact, the corresponding hypercohomology spectral sequence
with
$$
E^1_{ij}=H^j(\cF^{-i}),\quad i \leq 0,\quad j \geq 0
$$
abuts to zero.
Therefore, it suffices to prove that
$H^j(\cF^i)$ lives in degrees less than $i-j+p$ for $0<j<i$.

By the K\"unneth formula and Lemma~\ref{cohomologyvanish},
it suffices to prove (using Lemma~\ref{decompose}) that
$H^j(S_\lambda(n-i))=0$ for $n\ge i-j+p$, $0<j<i$.
Here $\lambda$ is any Young diagram with at most $p$ boxes in each row.
This follows from Proposition~\ref{main}.
\end{proof}

\begin{Corollary}\label{cor1}
$\Phi_*(D)$ is $(p+1)$-linear over $\cL$.
\end{Corollary}

\begin{proof}
We check that the complex \eqref{BigComplex} satisfies conditions
of Lemma~\ref{funny}. This follows from Proposition~\ref{LinearityD''}
and Lemma~\ref{linearityD'} (using Lemma~\ref{Segre}),
and from Proposition~\ref{BBW}.
\end{proof}

This concludes the proof of the Theorem \ref{inductiveversion}
\end{proof}

\section{Quadrics}\label{segrequadrics}
We make use of
Kapranov's Hilbert scheme realisation
of $\oM_{S^s}$ as the subscheme of
$\Hilb(\bP^{|S|-3})$ of Veronese curves (i.e.~stable $S^s$-pointed
rational curves embedded by global sections of $\omega_C(S^s)$) through
$|S^s|$ fixed general points in $\bP^{|S|-3}$, with
$$
\oM_S \subset \oM_{S^s} \times \bP^{|S| -3}
$$
the universal family. See \cite[Theorem 0.1]{Kapranov93a}.

We let $\pi = \pi_{S^s}$.

We consider the bundle $I_2$ on $\oM_{S^s}$ whose
fibre at $[C]$ is the vector space of quadrics
in $\bP^{n-3}$ vanishing on $C$. I.e.
$$
0 \rightarrow I_2 \rightarrow \sym_2(\bV) \rightarrow
\pi_*(\cO_{\oM_S} \otimes \psi_s^{\otimes 2}) \rightarrow 0
$$
where $\bV$ is the trivial bundle
$H^0(\oM_S,\psi_s)$. Note $I_2 = \cM^1$ of Lemma~\ref{defofMq}.

\begin{Definition}[Segre Quadrics]\label{4.0 definition}
Let $V,W$ be vector spaces,
$X,Y \in V$, $\sigma,\gamma \in W$
$$
Q(X,Y \in V,\sigma,\gamma \in W) :=
(X \otimes \sigma) \otimes (Y \otimes \gamma) -
(X \otimes \gamma) \otimes (Y \otimes \sigma) \in
(V \otimes W)^{\otimes 2}
$$
We will abuse notation and use the same symbol for the
image of $Q$ in $\Sym_2(V \otimes W)$.
\end{Definition}
\begin{Remark} It is well known that the homogenous ideal
of the Segre embedding
$$
\bP(V) \times \bP(W) \subset \bP(V \otimes W)
$$
is generated by Segre quadrics. \end{Remark}

The following is obvious from the definitions:

\begin{Lemma}\label{4.0}
For vector spaces $V,W,Z$ and
elements $X,Y \in V$, $\sigma,\gamma \in W$, $a,b \in Z$
$$
\begin{CD}
Q(X,Y \in V,\sigma \otimes a, \gamma \otimes b \in W \otimes Z) +
Q(X,Y \in V,\sigma \otimes b,\gamma \otimes a \in W \otimes Z) = \\
Q(X,Y \in V, \sigma, \gamma \in W) \otimes (a \otimes b + b \otimes a)
\end{CD}
$$
in
$$
(X \otimes (Y \otimes Z))^{\otimes 2} = (X \otimes Y)^{\otimes 2}\otimes
Z^{\otimes 2}.
$$
\end{Lemma}

\begin{Corollary}\label{4.0.1}
Let
$$
G' \subset \Sym_2(V \otimes W),\quad G \subset \Sym_2(V \otimes (W \otimes Z))
$$
be the subspaces generated by Segre quadrics. Then
$G ' \otimes \Sym_2(Z)$ is contained in the image of $G$
under the natural map
$$
\Sym_2(V \otimes (W \otimes Z)) \rightarrow \Sym_2(V \otimes W) \otimes
\Sym_2(Z).
$$
\end{Corollary}

\begin{Lemma}\label{lfs}
For any $F\subset S^s$,
$$H^0(\oM_{S^s},L_F) \otimes H^0(\oM_S,\psi_s) = H^0(\oM_S,L_F).$$
\end{Lemma}
\begin{proof} Follows from Lemma~\ref{residues} and  Lemma~\ref{old2.5}.
\end{proof}

\begin{Theorem}\label{quadrics}
Segre quadrics
$$
Q(X,Y \in H^0(\kappa_F),\sigma,\gamma \in H^0(L_F))
$$
for subsets $F \subset S$ generate the homogenous ideal of
$$
\oM_{S} \subset \bP(H^0(\kappa)^*).
$$
\end{Theorem}

\begin{proof}
We write $\cS^k$ for $\Sym_k$.

By Theorem~\ref{propertiesofkapparing} and Theorem~\ref{Koszulisquadratic}
the homogeneous ideal is generated by quadrics
so its enough to show Segre quadrics generate the
kernel of
\begin{equation}\label{4.1}
\cS^2(H^0(\kappa)) \rightarrow H^0(\kappa^{\otimes 2}) \cooltag
\end{equation}

\eqref{4.1} factors through
\begin{equation}\label{4.2}
 \cS^2(H^0(\kappa)) = \cS^2(H^0(\kappa_{S^s}) \otimes H^0(\psi_s)) \rightarrow
\cS^2(H^0(\kappa_{S^s})) \otimes \cS^2(H^0(\psi_s)) \cooltag
\end{equation}

Simple linear algebra shows the kernel of \eqref{4.2} is generated by
Segre quadrics, so its enough to show that images of Segre quadrics
generate the
kernel of
\begin{equation}\label{4.3}
\cS^2(H^0(\kappa_{S^s})) \otimes \cS^2(H^0(\psi_s)) \rightarrow
H^0(\kappa_{S^s}^{\otimes 2} \otimes \psi_s^{\otimes 2}) \cooltag
\end{equation}

\eqref{4.3} factors through
\begin{equation}\label{4.4}
\cS^2(H^0(\kappa_{S^s})) \otimes \cS^2(H^0(\psi_s)) \rightarrow H^0(\kappa_{S^s}^{\otimes 2}) \otimes
\cS^2(H^0(\psi_s)) \cooltag
\end{equation}

By induction, Corollary \ref{4.0.1}, and Lemma ~\ref{lfs} the
images of Segre quadrics  generate the
kernel of \eqref{4.4}, thus
it is enough to show that images of Segre quadrics generate the
kernel of
\begin{equation}\label{4.5}
H^0(\kappa_{S^s}^{\otimes 2}) \otimes \cS^2(H^0(\psi_s)) \rightarrow
H^0(\kappa^{\otimes 2}) \cooltag
\end{equation}

\eqref{4.5} is induced by pushforward and taking global sections from
\begin{equation}\label{4.6}
0 \rightarrow I_{\oM_S} \otimes \kappa_{S^s}^{\otimes 2} \otimes \psi_s^{\otimes 2}
\rightarrow \cO_{\oM_{S^s} \times \bP(H^0(\psi_s)^*)} \otimes
\kappa_{S^s}^{\otimes 2} \otimes \psi_s^{\otimes 2}  \rightarrow
\kappa_{\oM_S}^{\otimes 2} \rightarrow 0 \cooltag
\end{equation}

Pushing forward by $\pi_{S^s}$ we obtain

\begin{equation}\label{4.6a}
0 \rightarrow I_2 \otimes \kappa_{S^s}^{\otimes 2}
\rightarrow \kappa_{S^s}^{\otimes 2} \otimes \cS^2(\bV) \rightarrow
{\pi_{S^s}}_*(\kappa_{\oM_{S}}^{\otimes 2}) \rightarrow 0\cooltag
\end{equation}
where the right hand zero (which we will not use) is implied by
Corollary~\ref{npforVeronese}.

Now to prove the Theorem it is enough to show that
$$
H^0(I_2 \otimes \kappa_{S^s}^{\otimes 2}) \subset H^0(\kappa_{S^s}^{\otimes 2})
\otimes \cS^2(H^0(\psi_s))
$$
is generated by images of Segre quadrics. We prove this by induction
on $|S|$. Suppose first $|S| \geq 6$.

For $t \in S^s$ consider the commutative diagram.
$$
\begin{CD}
\oM_S @>>> \oM_{S^t} \\
@VVV     @VVV \\
\oM_{S^s} @>>> \oM_{S^{s,t}}
\end{CD}
$$
There is a natural inclusion
$\pi_{S^{s,t}}^*(I_2) \subset I_2$ (indeed a subbundle) which
induces a natural inclusion
$$
H^0(\oM_{S^{s,t}},I_2 \otimes \kappa_{S^{s,t}}^{\otimes 2})
\subset H^0(\oM_{S^s},I_2 \otimes \kappa_{S^{s,t}}^{\otimes 2})
$$
Which then gives a natural map

\begin{equation}\label{4.10}
H^0(\oM_{S^{s,t}},I_2 \otimes \kappa_{S^{s,t}}^{\otimes 2}) \otimes
\cS^2(H^0(\oM_{S^s},\psi_t))
\rightarrow H^0(\oM_{S^s},I_2 \otimes \kappa_{S^s}^{\otimes 2}). \cooltag
\end{equation}

Consider the following diagram:
$$
\begin{CD}
(H^0(\kappa_{S^t}) \otimes H^0(\pi_{S^s}^*(\psi_t)))^{\otimes 2}
@>e >> H^0(\kappa_S)^{\otimes 2} \\
 @| @| \\
H^0(\oM_{S^t},\kappa_{S^{s,t}})^{\otimes 2}  \otimes H^0(\oM_{S^s},\psi_t)^{\otimes 2}
\otimes H^0(\oM_{S^t},\psi_s))^{\otimes 2}
@. H^0(\kappa_{S^s})^{\otimes 2} \otimes H^0(\psi_s)^{\otimes 2} \\
@VVV @VVV \\
H^0(\oM_{S^s},\kappa_{S^{s,t}}^{\otimes 2}) \otimes \cS^2(H^0(\oM_{S^s},\psi_t))
\otimes \cS^2(H^0(\oM_{S^t},\psi_s))
@>f>> H^0(\kappa_{S^s}^{\otimes 2})
\otimes \cS^2(H^0(\oM_S,\psi_s))
\end{CD}
$$
where the maps are as follows:

We have natural identifications and inclusions
$$
\begin{CD}
\pi_{A^a}^*(\psi_b) = \psi_b(-\delta_{a,b}), ~ a \neq b \in A \\
H^0(\oM_{A^a},\psi_b) = H^0(\oM_A,\pi_{S^a}^*(\psi_b)) \subset H^0(\oM_A,\psi_b) \\
H^0(\oM_A,\kappa_{A^a}) \otimes H^0(\oM_A,\psi_a) = H^0(\oM_A,\kappa_A) \\
\cS^k(H^0(\oM_A,\psi_a)) = \cS^k(\bP^{|A|-3},\cO(1)) =
H^0(\oM_A,\psi_a^{\otimes k}) \\
H^0(\oM_S,\kappa_{S^{s,t}}^{\otimes k}) =
H^0(\oM_{S^s},\kappa_{S^{s,t}}^{\otimes k}) =
H^0(\oM_{S^t},\kappa_{S^{s,t}}^{\otimes k}).
\end{CD}
$$

$e$ is the composition
$$
(H^0(\kappa_{S^t}) \otimes H^0(\pi_{S^s}^*(\psi_t)))^{\otimes 2}
\subset
(H^0(\kappa_{S^t}) \otimes H^0(\oM_S,\psi_t))^{\otimes 2}
= (H^0(\oM_S,\kappa))^{\otimes 2}
$$
and the other maps are given in the obvious way by multiplication of
sections. One checks immediately that the diagram is commutative.

By Lemma \ref{4.13} and Lemma \ref{4.14} the images of maps
\eqref{4.10} over all $t \in S^s$ generate
$H^0(\oM_{S^s},I_2 \otimes \kappa_{S^s}^{\otimes 2})$.
By induction the image of
$$
H^0(\oM_{S^{s,t}},I_2 \otimes \kappa_{S^{s,t}}^{\otimes 2}) \subset
H^0(\kappa_{S^{s,t}}^{\otimes 2}) \otimes \cS^2(H^0(\oM_{S^t},\psi_s))
$$
is generated by images of Segre quadrics under the natural map
$$
H^0(\oM_{S^t},\kappa_{S^t})^{\otimes 2} =
(H^0(\oM_{S^{s,t}},\kappa_{S^{s,t}}) \otimes H^0(\oM_{S^t},\psi_s))^{\otimes 2}
\to H^0(\kappa_{S^{s,t}}^{\otimes 2}) \otimes \cS^2(H^0(\oM_{S^t},\psi_s)).
$$
Thus
$$
H^0(\oM_{S^{s}},I_2 \otimes \kappa_{S^{s}}^{\otimes 2}) \subset
H^0(\kappa_{S^{s}}^{\otimes 2}) \otimes \cS^2(H^0(\oM_{S},\psi_s))
$$
is generated by elements of form
$$
f(\oQ \otimes
\overline{(a \otimes b + b \otimes a}))
$$
where
$$
\oQ \in H^0(\kappa_{S^{s,t}}^{\otimes 2}) \otimes \cS^2(H^0(\oM_{S^t},\psi_s))
$$
is the image of a Segre quadric
$Q \in H^0(\kappa_{S^t})^{\otimes 2}$ and
$$
\overline{(a \otimes b + b \otimes a)}
\in \cS^2(H^0(\oM_{S^s},\psi_t))
$$
is the image of the symmetric tensor
$a \otimes b + b \otimes a \in H^0(\oM_{S^s},\psi_t)^{\otimes 2}$.
Now by the commutativity of the above diagram, Lemma~\ref{4.0},
and Lemma~\ref{lfs},
$f(\oQ \otimes \overline{(a \otimes b + b \otimes a)})$ is in the
span of images of Segre quadrics. This completes the induction step.

So now suppose $|S| =5$.
In this case $H^0(I_2 \otimes \kappa_{S^s}^{\otimes 2})$ is
two dimensional.
So its enough to show there are two Segre quadrics
whose images in $H^0(\kappa_{S^s}^{\otimes 2}) \otimes \cS^2(H^0(\psi_s))$ are
linearly independent.

Let $S =\{x,y,z,g,s\}$.
Let $R = \{x,y,s,g\}$. We will find a quadric
$$
Q = Q(X,Y \in H^0(\kappa_R),\sigma,\gamma \in H^0(\psi_z)))
$$
whose restriction
to the fibre $[C] = \delta_{x,y}\in  \oM_{S^s}$ is non-trivial,
but whose restriction to the fibre $[E] = \delta_{x,z}$ is
identically zero. The result will then follow by symmetry.

By the relations
$$
\begin{CD}
\kappa_R = \cO(\delta_{x,s} + \delta_{y,g})
= \cO(\delta_{y,s} + \delta_{x,g})  \\
\psi_z = \pi_{\{s,y,z,g\}}^*(\psi_z) + \cO(\delta_{x,z})
= \cO(\delta_{z,s} + \delta_{y,g} + \delta_{x,z}) =
\cO(\delta_{g,s} + \delta_{y,z} + \delta_{x,z} )
\end{CD}
$$
we may choose sections $X,Y \in H^0(\kappa_R)$,
$\sigma,\gamma \in H^0(\psi_z)$ with zero divisors
$$
\begin{CD}
Z(X)  = \delta_{x,s} + \delta_{y,g},\quad Z(Y) = \delta_{y,s} + \delta_{x,g}  \\
Z(\sigma) =  \delta_{z,s} + \delta_{y,g} + \delta_{x,z}, \quad
Z(\gamma) = \delta_{g,s} + \delta_{y,z} + \delta_{x,z}.
\end{CD}
$$
The restrictions of the sections
$$
X \otimes \sigma, Y \otimes \gamma, X \otimes \gamma
$$
to $H^0(C,\omega_C(S^s))$ have zero schemes
$x + z$, $y + g$, $x + g$, respectively. In particular
the restrictions of the three sections are linearly independent
and so they give a basis of (the three dimensional vector
space) $H^0(C,\omega_C(S^s))$. Now its clear that the
restriction of the quadric $Q$ to $\cS^2(H^0(C,\omega_C(S^s))$ is non-trivial.
Observe that the restriction of $Q$ to
$$
\cS^2(H^0(E,\omega_E(S^s))
$$
is identically zero.
Indeed the section $\sigma$ itself vanishes identically along $E$.
\end{proof}

\begin{Lemma}\label{temp1}
The map
\begin{equation}
\boplus_{a \in T^{s,t}} H^0(\pi_{T^a}^*(\cO(-\delta_{s,t})\otimes \kappa^{\otimes 2})
\to H^0(\cO(-\delta_{s,t})\otimes \kappa^{\otimes 2})
\end{equation}
is surjective for distinct $s,t \in T$, $|T| \geq 5$.
\end{Lemma}
\begin{proof}
Choose $x,y \in T^{s,t}$. Since $\delta_{x,s,t}$ and $\delta_{y,s,t}$
are disjoint, we have an exact sequence:
$$
0 \to \cO(-\delta_{s,t} - \delta_{s,t,x} - \delta_{s,t,y})
\to \cO(\pi_{S^x}^*(-\delta_{s,t})) \oplus \cO(\pi_{S^y}^*(-\delta_{s,t}))
\to \cO(-\delta_{s,t}) \to 0.
$$
Let $E := \delta_{s,t} + \delta_{s,t,x} + \delta_{s,t,y}$. Observe
$E$ satisfies the conditions of Lemma~\ref{boundarylemma}
(for $s \equiv t \not \equiv a$, $a \in T^{s,t,x,y}$). So the
sequence remains exact after tensoring by $\kappa^{\otimes 2}$ and
taking global sections. \end{proof}

\begin{Lemma}\label{temp2} The map
\begin{equation}
\boplus_{x \in T^s} H^0(\pi_{T^x}^*(V_{\psi_s}) \otimes \kappa^{\otimes 2})
\to H^0(V_{\psi_s} \otimes \kappa^{\otimes 2})
\end{equation}
is surjective, for $s \in T$, $|T| \geq 5$. \end{Lemma}
\begin{proof}
Choose $t \in T^{s}$.
We have by Lemma~\ref{firstVpsi} and Lemma~\ref{cohomologyvanish} a
commutative diagram with short exact rows (where we have omitted
the left and right zeros for reasons of space):
$$
\begin{CD}
\boplus_{x \in T^s} H^0(\pi_{T^{x,t}}^*(V_{\psi_s}) \otimes \kappa^{\otimes 2})
@>>> \boplus_{x \in T^s} H^0(\pi_{T^x}^*(V_{\psi_s}) \otimes \kappa^{\otimes 2})
@>>> \boplus_{x \in T^{s,t}} H^0(\cO(\pi_{T^x}^*(-\delta_{s,t})) \otimes \kappa^{\otimes 2})
\\ @VVV @VVV @VVV  \\
H^0(\pi_{T^t}^*(V_{\psi_s}) \otimes \kappa^{\otimes 2})
@>>> H^0(V_{\psi_s} \otimes \kappa^{\otimes 2}) @>>> H^0(\cO(-\delta_{s,t}) \otimes
\kappa^{\otimes 2})
\end{CD}
$$
Note the first column is surjective, since the sum on the upper left
includes the case $x = t$, on which the map is the identity. The right
column is surjective by the previous lemma. Thus the center column is
surjective. \end{proof}

\begin{Lemma}\label{4.13} Assume $|T| \geq 5$.
The map
$$
\boplus_{x \in T} H^0(\pi_{T^x}^*(I_2) \otimes \kappa^{\otimes 2})
\to H^0(I_2 \otimes \kappa^{\otimes 2})
$$
is surjective.
\end{Lemma}

\begin{proof}
Note that $I_2 = \cM^1$ of \eqref{defofMq}.
Choose $s \in T$.
We have by Theorem~\ref{b2.1} and Lemma~\ref{cohomologyvanish} a
commutative diagram with short exact rows (with left and right
zeros omitted):
$$
\begin{CD}
\boplus_{x \in T} H^0(\pi_{T^{x,s}}^*(I_2) \otimes \kappa^{\otimes 2})
@>>> \boplus_{x \in T} H^0(\pi_{T^x}^*(I_2) \otimes \kappa^{\otimes 2})
@>>> \boplus_{x \in T^{s}} H^0(V_{\psi_s} \otimes \kappa^{\otimes 2})
\\
@VVV @VVV @VVV  \\
H^0(\pi_{T^s}^*(I_2) \otimes \kappa^{\otimes 2})
@>>> H^0(I_2 \otimes \kappa^{\otimes 2}) @>>> H^0(V_{\psi_s} \otimes
\kappa^{\otimes 2}).
\end{CD}
$$
Note the first column is surjective since the upper left term
includes the case $x =s$ on which the map is the identity. The
right column is surjective by the previous lemma. The result follows. \end{proof}

\begin{Lemma}\label{4.14} For $t \in T$, $|T \geq 5|$ the
map
$$
H^0(\pi_{S^t}^*(I_2) \otimes \kappa_{S^t}^{\otimes 2})
\otimes \Sym_2(H^0(\psi_t)) \to H^0(\pi_{S^t}^*(I_2) \otimes \kappa^{\otimes 2})
$$
is surjective.
\end{Lemma}
\begin{proof}
Since $(\pi_{S^t})_*(\psi_t)$ is trivial, we have by the
projection formula
$$
H^0(\pi_{S^t}^*(F) \otimes \psi_t) = H^0(\pi_{S^t}^*(F)) \otimes H^0(\psi_t)
$$
for any vector bundle $F$ on $\oM_{T^t}$. So its enough to show
$$
H^0(\pi_{S^t}^*(I_2) \otimes \kappa_{S^t} \otimes \kappa)
\otimes H^0(\psi_t)) \to H^0(\pi_{S^t}^*(I_2) \otimes \kappa^{\otimes 2})
$$
is surjective. Note for any vector bundle $W$
$$
H^0(W) \otimes H^0(\psi_t) \to H^0(W \otimes \psi_t)
$$
is surjective so long as $H^1(W \otimes V_{\psi_t}) =0$.
$I_2 \otimes \kappa$ is extended from globally generated
line bundles by Theorem~\ref{b2.1} and Lemma~\ref{mainfiltration}.
So we have the necessary vanishing by Corollary~\ref{mainvanishing}.
\end{proof}

\end{document}